\documentclass{jnmp03c}

\usepackage{amstext}

\begin{document}
\setcounter{page}{59}
\renewcommand{\evenhead}{A Sergeev}
\renewcommand{\oddhead}{Superanalogs of the Calogero Operators
and Jack  Polynomials}

\thispagestyle{empty}

\FistPageHead{1}{\pageref{sergeev-firstpage}--\pageref{sergeev-lastpage}}{Letter}

\copyrightnote{2001}{A Sergeev}

\Name{Superanalogs of the Calogero Operators\\
 and Jack  Polynomials}\label{sergeev-firstpage}

\Author{A SERGEEV}

\Adress{Department of Mathematics, University of
Stockholm\\
Roslagsv.  101, Kr\"aftriket hus 6, S-106 91, Stockholm, Sweden\\
E-mail: mleites@matematik.su.se\\[2mm]
 On leave of absence from\\
 Balakovo Institute of Technique of Technology and Control Branch\\
 Saratov Technical University,  Balakovo, Saratov Region, Russia}

\Date{Received November 12, 2000; Accepted December 12, 2000}

\begin{abstract}
\noindent
 A depending on a complex parameter $k$ superanalog
${\mathcal S}{\mathcal L}$ of Calogero operator is constructed; it is related with the
root system of the Lie superalgebra ${\mathfrak{gl}}(n|m)$. For $m=0$ we obtain
the usual Calogero operator; for $m=1$ we obtain, up to a change of
indeterminates and parameter $k$ the operator constructed by Veselov,
Chalykh and Feigin~[2,~3].
For $k=1, \frac12$ the operator ${\mathcal S}{\mathcal L}$ is the
radial part of the 2nd order Laplace operator for the symmetric
superspaces corresponding to pairs $(GL(V)\times GL(V), GL(V))$ and
$(GL(V), OSp(V))$, respectively. We will show that for the generic $m$
and $n$ the superanalogs of the Jack polynomials constructed by Kerov,
Okunkov and Olshanskii~[5]
are eigenfunctions of ${\mathcal S}{\mathcal L}$; for $k=1, \frac12$
they coinside with the spherical functions corresponding to the above
mentioned symmetric superspaces.
We also study the inner product induced by Berezin's integral on these
superspaces.
\end{abstract}

\noindent The Hamiltonian of the quantum Calogero problem is of
the form
\be
{\mathcal L}=\sum\limits_{i=1}^n \left(\frac{\p}{\p
t_{i}}\right)^2- \frac12 k(k-1) \sum\limits_{i<j}
\frac{\omega^2}{\sinh^2\frac{\omega}{2}(t_{i}-t_{i})} .
\ee
In
this form it is a particual case (corresponding to the root system
$R$ of ${\mathfrak{gl}}(n)$) of the operator constructed in the
famous paper by Olshanetsky and Perelomov~\cite{OP} \be {\mathcal
L}=\Delta -\mathop{\sum}\limits_{\alpha \in R^+}k_{\alpha}
(k_{\alpha}-1) \frac{(\alpha,
\alpha)}{\left(e^{\frac12\alpha}-e^{-\frac12\alpha}\right)^2} .
\ee Veselov, Feigin and Chalykh~\cite{CFV1} suggested the
following generalization of operator (1) 
\be \ba{l} \ds {\mathcal
L}'=\sum\limits_{i=1}^n \left(\frac{\p}{\ p t_{i}}\right)^2+
\left(\frac{\p}{\p t_{n+1}}\right)^2- \frac12 k(k-1)
\sum\limits_{i<j}
\frac{\omega^2}{\sinh^2\frac{\omega}{2}(t_{i}-t_{i})}
\vspace{3mm}\\
\ds \phantom{{\mathcal L}'=} + \frac12(k-1)\sum\limits_{i=1}^n
\frac{\omega^2}{\sinh^2\frac{\omega}{2}\left(t_{i}-\sqrt{-k}t_{n+1}\right)}.
\ea \ee 
It is known (\cite{LV}) that eigenfunctions of operator
(1) can be expressed in terms of Jack polnomials
$P_{\lambda}(x_{1},\ldots , x_{n}; k)$, where $\lambda$ is a
partition of $n$.  (For definition and properties of Jack
polynomials see \cite{M,St}.)  It is known~(\cite{M}) that for
$k=1, \frac12, 2$ (our $k$ is inverse of $\alpha$, the parameter
of Jack polynomials Macdonald uses in~\cite{M}) Jack polynomials
are interpreted as spherical functions on symmetric spaces
corresponding to pairs $(GL\times GL, GL)$, $(GL, SO)$ and $(GL,
Sp)$, respectively.  In these cases the corresponding operators
are radial parts of the corresponding second order Laplace
operators.

\medskip

\noindent
{\bf Superroots of $\pbf{{\mathfrak{gl}}(n|m)}$:}

Let $I=I_{{\bar 0}}\coprod I_{{\bar 1}}$
be the union of the ``even'' indices $I_{{\bar 0}}=\{1, \dots , n\}$ and ``odd''  indices
$I_{{\bar 1}}=\{\bar 1, \dots, \overline m\}$. Let $\dim V=(n|m)$ and $e_{1},
\dots, e_{n}, e_{\bar 1}, \dots , e_{\overline m}$ be a basis of $V$ such
that the parity of each vector is equal to that of its index. Let
 $\varepsilon_{1},
\dots, \varepsilon_{n}, \varepsilon_{\bar 1}, \dots ,
\varepsilon_{\bar  m}$ be the left dual basis of~$V^{*}$. Then the
set of roots can be described as follows: $R=R_{11}\coprod
R_{22}\coprod R_{12}\coprod R_{21}$, where \be \ba{l} \ds
R_{11}=\{\varepsilon_{i}-\varepsilon_{j}\mid i, j\in I_{{\bar
0}}\},\qquad R_{22}=\{\varepsilon_{i}-\varepsilon_{j}\mid i, j\in
I_{{\bar 1}}\},
\vspace{2mm}\\
R_{12}=\{\varepsilon_{i}-\varepsilon_{j}\mid i\in I_{{\bar 0}}, \; j\in I_{{\bar 1}}\},\qquad
R_{21}=\{\varepsilon_{i}-\varepsilon_{j}\mid i\in I_{{\bar 1}}, \; j\in I_{{\bar 0}}\}.
\ea
\ee
On $V^*$, define the depending on parameter $k$ inner product by
setting
\be
(v_{1}^*, v_{2}^*)_{k}=\mathop{\sum}\limits_{i=1}^n
v_{1}^*(e_{i})v_{2}^*(e_{i})-k\mathop{\sum}\limits_{j=1}^m
v_{1}^*(e_{\bar j})v_{2}^*(e_{\bar j})
\ee
and set $\rho_{(k)}=k\rho_{1}+\frac1k\rho_{2}-\rho_{12}$, where
\[
\rho_{1}=\frac12\sum\limits_{\alpha\in R_{11}^+}  \alpha; \qquad
\rho_{2}=\frac12\sum\limits_{\beta\in R_{22}^+}  \beta; \qquad
\rho_{12}=\frac12\sum\limits_{\gamma\in R_{12}}  \gamma.
\]

Set $\partial_{i}(e^{v^*})= v^*(e_i)e^{v^*}$, $\partial_{\bar
j}(e^{v^*})= v^*(e_{\bar j})e^{v^*}$.  Define the superanalog of the
Calogero operator to be
\be
\ba{l}
\ds {\mathcal S}{\mathcal L}=\sum\limits_{i=1}^n \partial_{i}^2-
k \sum\limits_{\bar j=1}^m \partial_{\bar j}^2
-k(k-1)\sum\limits_{\alpha\in R_{11}^+}\frac{(\alpha, \alpha)_{k}}
{\left(e^{\frac12\alpha}-e^{-\frac12\alpha}\right)^2}
\vspace{3mm}\\
\ds \phantom{{\mathcal S}{\mathcal L}=}
+  \frac{1}{k}\left(\frac{1}{k}-1\right)
\sum\limits_{\beta\in R_{22}^+}
\frac{(\beta, \beta)_{k}}
{\left(e^{\frac12\beta}-e^{-\frac12\beta}\right)^2}-
2\sum\limits_{\gamma\in R_{12}}
\frac{(\gamma, \gamma)_{k}}
{\left(e^{\frac12\gamma}-e^{-\frac12\gamma}\right)^2}.
\ea
\ee

In order to describe the eigenfunctions of ${\mathcal S}{\mathcal L}$,
it is convenient to
present it in terms of operator ${\mathcal M}$ described below.  Set
\be
\delta^{(k)}=\prod\limits_{\alpha\in R_{11}^+}
\left(e^{\frac12\alpha}-e^{-\frac12\alpha}\right)^k
\prod\limits_{\beta\in R_{22}^+}
\left(e^{\frac12\beta}-e^{-\frac12\beta}\right)^{1/k}
\prod\limits_{\gamma\in R_{12}}
\left(e^{\frac12\gamma}-e^{-\frac12\gamma}\right)^{-1}.
\ee
Set
\[
{\mathcal M}=\left(\delta^{(k)}\right)^{-1}\left ({\mathcal L}-(\rho_{(k)},
\rho_{(k)})_{k}\right )\delta^{(k)}.
\]

\noindent
{\bf Lemma.} {\it The explicit form of ${\mathcal M}$ is
\be
{\mathcal M}=\sum\limits_{i=1}^n \partial_{i}^2-
k \sum\limits_{\bar j=1}^m \partial_{\bar j}^2
+k \sum\limits_{\alpha\in R_{11}^+}\frac{e^{\alpha}+1}
{e^{\alpha}-1}\partial_{\alpha}-
 \sum\limits_{\beta\in R_{22}^+}
\frac{e^{\beta}+1}{e^{\beta}-1}\partial_{\beta}-
2 \sum\limits_{\gamma\in R_{12}}
\frac{e^{\gamma}+1}
{e^{\gamma}-1}\partial_{\gamma, k},
\ee
where
\be
\ba{l}
\ds  \partial_{\alpha}=\partial_{i}-\partial_{j} \qquad \text{for}
\quad \alpha=\varepsilon_{i}-\varepsilon_{j},
\vspace{2mm}\\
\ds \partial_{\beta}=\partial_{\bar i}-\partial_{\bar j} \qquad \text{for}
\quad \beta=\varepsilon_{\bar i}-\varepsilon_{\bar j},
\vspace{2mm}\\
\ds \partial_{\gamma, k}=\partial_{i}-\partial_{\bar j}\qquad \text{for}
\quad \gamma=\varepsilon_{i}-k\varepsilon_{\bar j}.
\ea
\ee
In terms of new indeterminates $x_{i}=e^{\varepsilon_{i}}$ and $y_{\bar
j}=e^{\varepsilon_{\bar j}}$ the operator ${\mathcal M}$ takes the form
\be
\ba{l}
\ds {\mathcal M}=\sum\limits_{i=1}^n \left (x_{i}\frac{\p}{\p x_{i}}\right )^2-
k  \sum\limits_{\bar j=1}^m \left (y_{j}\frac{\p}{\p y_{j}}\right )^2+
k \sum\limits_{1\leq i<j\leq n}\frac{x_{i}+x_{j}}
{x_{i}-x_{j}}\left (x_{i}\frac{\p}{\p x_{i}}-x_{j}\frac{\p}{\p x_{j}}\right )
\vspace{3mm}\\
\ds \phantom{{\mathcal M}=} -\!
\sum \limits_{1\leq i<j\leq n}\frac{y_{i}+y_{j}}
{y_{i}-y_{j}}\left (y_{i}\frac{\p}{\p y_{i}}-y_{j}\frac{\p }{\p y_{j}}\right )-\!\!
\sum \limits_{1\leq i\leq n, \; 1\leq j\leq m} \! \frac{x_{i}+y_{j}}
{x_{i}-y_{j}}\left (x_{i}\frac{\p}{\p x_{i}}-y_{j}\frac{\p }{\p y_{j}}\right ).
\ea\hspace{-15mm}
\ee}

Following Kerov, Okunkov and Olshanskii \cite{KOO}, determine
superanalogs of Jack polynomials. Let
\[
S^p=x_{1}^{p}+\dots +x_{n}^{p}\qquad \text{and} \qquad
S^\mu=S^{\mu_{1}}\cdots S^{\mu_{l}} \quad \text{for any partition
$\mu$ of $n$}
\]
and let $P_{\lambda}(x; k)=\sum \chi_{\mu}^\lambda S^\mu$ be the
decomposition of the classical Jack polynomials into sums of powers.
Set further
\[
S^{p, k}=\sum\limits_{1\leq i\leq
n}x_{i}^{p}-\frac1k \sum\limits_{1\leq j\leq m}y_{j}^{p}\quad \text{and}
\quad S^{\mu, k}=S^{\mu_{1}, k}\dots S^{\mu_{l}, k} \quad
\text{for any partition $\mu$ of $n$}.
\]
Then the superanalogs of Jack polynomials are of the form
\be
P_{\lambda}(x, y; k)=\sum \chi_{\mu}^\lambda(k) S^{\mu, k}.
\ee

\noindent {\bf Theorem 1.} {\it The polynomials $P_{\lambda}(x, y;
k)$ defined by eq.~$(9)$ are eigenfunctions of operator~$(10)$.}

\medskip

\noindent
{\bf Spherical functions:}

In this paper we adopt an algebraic
approach to the theory of spherical functions.

Let ${\mathfrak g}$ be a finite dimensional Lie superalgebra, $U({\mathfrak g})$ its
enveloping algebra, ${\mathfrak b}\subset {\mathfrak g}$ a subalgebra.  Let
$\pi:{\mathfrak g}\longrightarrow{\mathfrak{gl}}(V)$ be an irreducible representation and $V^*$ the
dual module.  If $v\in V$ is a nonzero ${\mathfrak b}$-invariant
vector, then there exists a nonzero vector $v^*\in V^*$ which is also
${\mathfrak b}$-invariant. The matrix coefficient
\[
\theta^\pi(v^*, v)\in U({\mathfrak g})^*
\]
will be called the {\it spherical function associated with the
triple} $(\pi, v^*, v)$.

Let $l\in U({\mathfrak g})^*$ be a left and right  ${\mathfrak b}$-invariant functional, i.e.,
\[
l(xuy)=l(u)\qquad \text{for any $x,y \in {\mathfrak b}$ and $u\in
U({\mathfrak g})$}.
\]
If $z\in Z({\mathfrak g})$, then $L^*(z)l$, where $L^*$ is left coregular
representation of ${\mathfrak g}$, is also a left and right
${\mathfrak b}$-invariant functional.

Let ${\mathfrak g}={\mathfrak{gl}}(V)\oplus {\mathfrak{gl}}(V)$ and
let ${\mathfrak b}_{1}\simeq {\mathfrak g}(V)$ be the first summand of
${\mathfrak g}$, whereas ${\mathfrak b}\simeq {\mathfrak g}(V)$ is the
diagonal subalgebra, i.e., ${\mathfrak b}=\{(x, x)\mid x\in {\mathfrak
g}(V)\}$.  Let ${\mathfrak h}$ be the Cartan subalgebra of
${\mathfrak{gl}}(V)$, let $\lambda$ be a partition of $l\in {\mathbb
N}$ and $V^\lambda$ an irreducible ${\mathfrak{gl}}(V)$-module in
$V^{\otimes l}$, corresponding to $\lambda$, see~\cite{S1}.

The ${\mathfrak g}$-module $W^\lambda=V^\lambda\otimes (V^\lambda)^*$ is
irreducible and contains a unique, up to a constant factor, invariant
vector $v_{\lambda}$.  The dual module $(W^\lambda)^*$ contains a
similar vector $v_{\lambda}^*$.  Let
$\varphi_{\lambda}=\theta^\pi(v_{\lambda}^*, v_{\lambda})$ be the
corresponding spherical function.

Let $\dim V=(n|m)$, and let $I_{{\bar 0}}=\{1, \dots , n\}$ and
$I_{{\bar 1}}=\{\bar 1, \dots, \overline m\}$; let $\{e_{ij}\mid i, j\in I=I_{{\bar 0}}
\coprod I_{{\bar 1}}\}$ be the basis of ${\mathfrak{gl}}(V)$ consisting of matrix units.
Set
\[
C_{2}=\sum\limits_{i\in
I_{{\bar 0}}}e_{ii}^{2}- \sum\limits_{j\in I_{{\bar 1}}}e_{jj}^{2}+
\sum\limits_{i, j\in I_{{\bar 0}};\; i\neq j}e_{ij}e_{ji}-
 \sum\limits_{i, j\in I_{{\bar 1}};\; i\neq j}e_{ij}e_{ji}-
\sum \limits_{i\in I_{{\bar 0}};\; j\in I_{{\bar 1}}}e_{ij}e_{ji}-
 \sum\limits_{i\in I_{{\bar 1}};\; j\in I_{{\bar 0}}}e_{ij}e_{ji}.
\]
As is easy to verify, $C_{2}$ is a central element in the enveloping
algebra of ${\mathfrak{gl}}(V)$
and ${\mathfrak g}$, if ${\mathfrak{gl}}(V)$ is embedded as the first
summand.

\bigskip

\noindent {\bf Theorem 2.} {\it
\begin{itemize}
\topsep0mm
\partopsep0mm
\parsep0mm
\itemsep0mm
\item[i)] Every left and right  invariant
functional $l\in U({\mathfrak g})^*$ is uniquely determined by its restriction
onto $S({\mathfrak h})\subset S({\mathfrak b}_{1})$.
\item[ii)] Let $S({\mathfrak h})^{inv}$ be the set of restrictions of left and right
invariant functional $l\in U({\mathfrak g})^*$
onto $S({\mathfrak h})\subset S({\mathfrak b}_{1})$.
Then for every $z\in Z({\mathfrak g})$ there exists a uniquely determined operator
$\Omega_{z}$ on $S({\mathfrak h})^{inv}$ (the radial part of $z$). It is
determined from the formula
\[
(\Omega_{z}l')(u)=(L^*(z)l)(u)\qquad \text{for any $l'\in
S({\mathfrak h})^{inv}$ and any its extension $l\in S({\mathfrak
g})$}.\vspace{-2mm}
\]
\item[iii)] The above defined operator $\Omega_{C_{2}}$ corresponding
to $C_{2}$ coinsides with the operator ${\mathcal M}$ determined by formula
$(10)$ for $k=1$.
\item[iv)] The functions $\varphi_{\lambda}$, as functionals on
$S({\mathfrak h})$, coinside, up to a constant factor, with Jack polynomials
$P_{\lambda}(x, y; 1)$, where $x_{i}=e^{\varepsilon _{i}}$ for $i\in I_{{\bar 0}}$ and
$y_{j}=e^{\varepsilon _{j}}$ for $j\in I_{{\bar 1}}$.
\end{itemize}}

\medskip

Let ${\mathfrak g}={\mathfrak{gl}}(V)$, where $\dim V=(n|m)$ and $m=2r$
is even. Let ${\mathfrak b}={\mathfrak{osp}}(V)$ be the
ortho\-symp\-lectic Lie subsuperalgebra in ${\mathfrak{gl}}(V)$
which preserves the tensor \be \sum\limits_{i\in I_{{\bar
0}}}e_{i}^{*}\otimes e_{i}^{*}+ \sum \limits_{j\in I_{{\bar
1}}}\left(e_{\bar j}^{*}\otimes e_{\overline{j+r}}^{*}+
e_{\overline{j+r}}^{*}\otimes e_{\bar j}^{*}\right). \ee Let
$\psi$ be an involutive automorphism of ${\mathfrak g}$ that
singles out ${\mathfrak{osp}}(V)$:
\[
{\mathfrak{osp}}(V)=\{x\in{\mathfrak{gl}}(V)\mid \psi(x)= -x\}.
\]
Let $V^\lambda$ be a ${\mathfrak g}$-module.  By~\cite{S2},
$V^\lambda$ contains a ${\mathfrak b}$-invariant vector $\tilde v_{\lambda}$ if
and only if $\lambda =2\mu$ and all its rows are of even length.  The
vector $\tilde v_{\lambda}^*\in (V^\lambda)^*$ is similarly defined.
Let $\tilde\varphi_{\lambda}=(v_{\lambda}^*, v_{\lambda})$ be the
corresponding matrix coefficient. Set
${\mathfrak h}^+=\{x\in {\mathfrak h}\mid \psi (x)=x\}$,
where ${\mathfrak h}\subset {\mathfrak g}$ is Cartan
subalgebra.

%\bigskip
\pagebreak

\noindent
{\bf Theorem 3.} {\it
\begin{itemize}
\topsep0mm
\partopsep0mm
\parsep0mm
\itemsep0mm
\item[i)] Every left and right  invariant functional on
$U({\mathfrak g})$ is uniquely determined by its restriction onto $S({\mathfrak h}^+)$.
\item[ii)] Let $S({\mathfrak h})^{inv}$ be the set of restrictions of left and right
invariant functionals.  Then for every $z\in Z({\mathfrak g})$ there exists a uniquely
determined operator $\Omega_{z}$ on $S({\mathfrak h})^{inv}$ (the radial part of
$z$).  It is determined from the formula
\[
(\Omega_{z}l')(u)=(L^*(z)l)(u)\qquad \text{for any $l'\in
S({\mathfrak h})^{inv}$ and any its extension $l\in S({\mathfrak
g})$}. \vspace{-2mm}
\]
\item[iii)] The operator $\Omega_{C_{2}}$ corresponding
to $C_{2}$ coinsides with the operator ${\mathcal M}$ determined by formula
$(10)$ for $m=r$ and $k=\frac12$.
\item[iv)] The functions $\tilde \varphi_{\lambda}$, as functionals on
$S({\mathfrak h})$, coinside, up to a constant factor, with Jack polynomials
$P_{\lambda}(x, y; \frac12)$, where $\lambda = 2\mu$, $x_{i}=e^{2\varepsilon
_{i}}$ for $1\leq i\leq n$ and $y_{j}=e^{2\varepsilon _{j}}$ for $1\leq j\leq
r$.
\end{itemize}}

\noindent
{\bf Invariant functional $F$ (the Berezin integral):}

For every ${\mathfrak g}$-module $W$, define in
$U({\mathfrak g})^*$ the subspace $C(W)$ consisting of the linear hull of the
matrix coefficients of $W$.  Denote by ${\mathfrak A}_{n, m}$ the subalgebra of
$U({\mathfrak g})^*$ generated by the matrix coefficients of the identity
representation $V$ of ${\mathfrak g}={\mathfrak{gl}}(V)$ and its dual.

\bigskip

\noindent
{\bf Theorem 4.}{\it
\begin{itemize}
\topsep0mm
\partopsep0mm
\parsep0mm
\itemsep0mm
\item[i)] On ${\mathfrak A}_{n, m}$, there exists a unique up to
a constant factor nontrivial left and right  invariant (with
respect to the left and right coregular representations) linear
functional~$F$.
\item[ii)] On ${\mathfrak A}_{n, m}$, define the inner product $\langle l_{1},
l_{2}\rangle =F(l_{1}^tl_{2})$, where $l\mapsto l^t$ is the principal
automorphism of $U({\mathfrak g})^*$. Then $\langle l_{1},
l_{2}\rangle =0$ for any $l_{1}\in C(V^\lambda)$,
$l_{2}\in C(V^\mu)$ and $\lambda\neq \mu$.
\item[iii)] If $\dim V^\lambda_{{\bar 0}}\neq \dim V^\lambda_{{\bar 1}}$, then
$\langle l_{1}, l_{2}\rangle =0$ for any $l_{1}, l_2\in C(V^\lambda)$.
\end{itemize}}

\medskip

{\bf Acknowledgements.}
I am thankful to D~Leites and G~Olshanskii for encouragement
and help.

\label{sergeev-lastpage}


\begin{thebibliography}{99}
\small
\topsep0mm
\partopsep0mm
\parsep0mm
\itemsep0mm

\bibitem{Di}
Dixmier J, Alg\` ebres Envellopentes, Gautier-Villars, Paris,
1974; Enveloping algebras, AMS, 1996.
\bibitem{CFV1}
Veselov A, Feigin M and  Chalykh O,  New Integrable Deformations of the
Quantum Calogero--Moser Problem, {\it Uspehi Mat. Nauk}, 1996, V.51, N~3,
185--186 (in Russian).
\bibitem{CFV2}
Chalykh O, Feigin M and Veselov A, Multidimesional Baker --Akhiezer
Functions and Huygens' Principle, {\it Comm. Math. Phys.}, 1999, V.206, N~3, 533--566.

\bibitem{K}
Kirillov A (jr.), Traces of Intertwining Operators and Macdonald's
Polynomials, Ph.D. thesis, q-alg/9503012.
\bibitem{KOO}
Kerov S, Okunkov A and Olshanskii G, The Boundary of the Young Graph
with Jack Edge Multipliers, {\it Internat. Research Notes}, 1998, N~4,
173--199.
\bibitem{LV}
Lapoint L and Vinet L, Exact Operator Solutions of the
Calogero--Sutherland Model, {\it Comm. Math. Phys.}, 1996, V.178, N~2,
425--452.
\bibitem{M}
Macdonald I, Symmetric Functions and Hall Polynomials, 2nd
Edition, Oxford Univ. Press, 1995.

\bibitem{OP}
Olshanetsky M and  Perelomov A, Quantum Integrable Systems Related to
Lie Algebras, {\it Phys. Rep.}, 1983, V.94, N~6, 313--404.
\bibitem{S1}
Sergeev A N, The Tensor Algebra of the Identity Representation as a
Module over the Lie Superlgebras $GL(n,m)$ and $Q(n)$, {\it Math. USSR
Sbornik}, 1985, V.51, 419--427.
\bibitem{S2}
Sergeev A, An Analog of the Classical Invariant Theory for Lie
Superalgebras. I, II, math.RT/9810111; math.RT/9904079; {\it Michigan J.
Mathematics}, 2001 (to appear).

\bibitem{St}
Stanley R, Some Combinatorial Properties of Jack Symmetric
Functions, {\it  Adv. Math.}, 1996, V.77, 76--115.
\end{thebibliography}
\end{document}